\documentclass{article}

\newcounter{note}
\newenvironment{note}[1][\hspace{-1.0ex}]%
 {\par\addvspace{2mm}\noindent\refstepcounter{note}\textbf{Remark~\thenote\hspace{1.0ex}{\rm#1}.~}\rm}%
 {\par\addvspace{2mm}\rm}
\newcounter{theorem}
\newenvironment{theorem}[1][\hspace{-1.0ex}]%
 {\par\addvspace{2mm}\noindent\refstepcounter{theorem}\textbf{Theorem~\thetheorem\hspace{1.0ex}{\rm#1}.~}\sl}%
 {\par\addvspace{2mm}\rm}
\newcounter{lemma}
\newenvironment{lemma}[1][\hspace{-1.0ex}]%
 {\par\addvspace{2mm}\noindent\refstepcounter{lemma}\textbf{Lemma~\thelemma\hspace{1.0ex}{\rm#1}.~}\sl}%
 {\par\addvspace{2mm}\rm}
\newcounter{prop}
\newenvironment{prop}[1][\hspace{-1.0ex}]%
 {\par\addvspace{2mm}\noindent\refstepcounter{prop}\textbf{Proposition~\theprop\hspace{1.0ex}{\rm#1}.~}\sl}%
 {\par\addvspace{2mm}\rm}
\newcounter{corol}
 {\par\addvspace{2mm}\noindent\refstepcounter{corol}\textbf{Corollary~\thecorol\hspace{1.0ex}{\rm#1}.~}\sl}%
 {\par\addvspace{2mm}\rm}

\newcommand\proofr[1][\hspace{-1em}]{\par{\em Proof\hspace{1em}#1:~}}
\providecommand\proofend{$\bigtriangleup$\vspace{0.3em}\par}

\title{On the binary codes with parameters of doubly-shortened $1$-perfect codes}
\author{Denis S. Krotov}
\begin{document}
\maketitle
\begin{abstract}
We show that any binary $(n=2^m-3, 2^{n-m}, 3)$ code $C_1$ is a part
of an equitable partition (perfect coloring) $\{C_1,C_2,C_3,C_4\}$ of the $n$-cube
with the parameters $((0,1,n-1,0) (1,0,n-1,0) (1,1,n-4,2) (0,0,n-1,1))$.
Now the possibility to lengthen the code $C_1$ to a $1$-perfect code
 of length $n+2$ is equivalent to the
possibility to split the part $C_4$ into two distance-$3$ codes or, equivalently,
to the biparticity of the graph of distances $1$ and $2$ of $C_4$.
In any case, $C_1$ is uniquely embeddable in a twofold $1$-perfect code of length
$n+2$ with some structural restrictions, where by a twofold $1$-perfect code
we mean that any vertex of the space is within radius $1$ from exactly two
codewords.
\end{abstract}

%\section{Introduction}
The hypercube $H^n=(V(H^n),E(H^n))$ of dimension $n$ is the graph
whose vertices are the all binary $n$-words,
two words being adjacent if and only if
they differ in exactly one position.

$d(\cdot,\cdot)$ -- the Hamming distance, i.e., the natural graph distance in $H^n$.

$\bar 0 = 0\ldots 0$ (the all-zero word), $\bar 1 = 1\ldots 1$
(the all-one word).

A binary code $C$ of length $n$ and code (or minimal) distance $d$, or $(n,|C|,d)$ code,
is a subset of $V(H^n)$ such that $d(\bar x, \bar y) \geq d$ for any different
$\bar x$ and $\bar y$ from $C$.

A partition $\{C_1,\ldots,C_r\}$ of $V(H^n)$ into $r$ nonempty parts
is said to be \emph{equitable} with parameters $(S_{ij})_{i,j=1}^n$
if for every $i,j\in\{1,\ldots,r\}$
every vertex $\bar x$ from $C_i$ has exactly
$S_{ij}$ neighbors from $C_j$
(the corresponding $r$-valued function on $V(H^n)$
is known as a \emph{perfect coloring}).

A binary code $C\subset V(H^n)$ is said to be \emph{$1$-perfect}
if every vertex $\bar x\in  V(H^n)$ is
at the distance $0$ or $1$ from exactly one codeword. Equivalently,
$\{C,V(H^n)\setminus C\}$ is an equitable partition with parameters
$((0,n)(1,n-1))$. Equivalently, $C$ is a $(2^m-1,2^{2^m-m-1},3)$ code, $n=2^m-1$.

We will say that a multiset $B\subset V(H^n)$ is a \emph{twofold $1$-perfect code}
if every vertex $\bar x\in  V(H^n)$ is
at the distance $0$ or $1$ from exactly two codewords of $B$.
We will say that a multiset $B\subset V(H^n)$ is \emph{splittable}
if it can be represented as the (multiset)
union of two distance-$3$ codes; otherwise $B$ is \emph{unsplittable}.
The existence of unsplittable twofold $1$-perfect codes was proved in \cite{KroPot:nonsplittable}.

We say that a code $C'$ if obtained by \emph{shortening} from a code $C\subset V(H^n)$
if $C' = \{\bar x\in V(H^{n-1}) \mid \bar x 0 \in C\}$. Respectively,
$C''$ is \emph{doubly-shortened} from $C$ if
$C'' = \{\bar x\in V(H^{n-2}) \mid \bar x 00 \in C\}$. (Here and elsewhere, for $\bar x = x_1x_2...x_n$, by $\bar x0$ we mean the concatenation of $\bar x$
with $0$, i.e., the word $x_1x_2...x_n0$;
similarly we define $\bar x1$, $\bar x00$, $\bar x01$, \ldots; we also expand this notation for
sets of words, e.g., $C0 = \{\bar x 0 \mid \bar x \in C\}$.)

It is known \cite{BesBro77} that shortened and doubly-shortened
(and even triply-shortened) $1$-perfect codes have the maximal
cardinality among all the codes of the same length and code distance $3$.
The question \cite{EV:98} is: can every code with such parameters
( $(2^m-2,2^{2^m-m-2},3)$ or $(2^m-3,2^{2^m-m-3},3)$ ) be represented
as a shortened or doubly-shortened $1$-perfect code?

For $(2^m-2,2^{2^m-m-2},3)$ codes the question is solved \cite{Bla99}.
In fact, such a code $C_1$ generates an equitable partition
$\{C_1, C_2, C_3\}$ with parameters $((0,n,0)(1,n-2,1)(0,n,0))$.
Then, the code
$$C = C_1 0 \cup  C_3 1$$
is $1$-perfect.

In this paper we prove that a $(2^m-3,2^{2^m-m-3},3)$ code $C_1$
generates an equitable partition
$\{C_1, C_2, C_3, C_4\}$ with parameters
$((0,1,n-1,0) (1,0,n-1,0) (1,1,n-4,2) (0,0,n-1,1))$.
If the code $C_4$ is splittable into two distance-$3$ codes $C'$ and $C''$,
then the code
$$C =  C_1 00 \cup  C_2 11 \cup
  C'01 \cup  C'' 10
$$
is $1$-perfect. However, the problem of splittability of $C_4$ remains open.
So, the problem of embedding $C_100$ in a $1$-perfect code is unsolved;
although, $C_100$ is proved to be embedded in twofold $1$-perfect codes
$$2\times C_1 00 \cup 2\times C_2 11 \cup C_401 \cup C_4 10$$
and
$$ C_1 00  \cup  C_2 00 \cup C_1 11  \cup  C_2 11 \cup C_401 \cup C_4 10$$
(Theorems~\ref{th:twofold2} and~\ref{th:twofold1}),
whose splittability is equivalent to the splittability of $C_4$.

\section{Notation and basic facts}
Let $C_1$ be a binary code of length $n=2^m-3$, cardinality
$2^{n-m}$, and minimal distance $3$.

Denote
\begin{eqnarray} \label{eq:C2}
C_2 &=& C_1 + \bar 1 = \{ \bar x \mid \bar x + \bar 1 \in C_1 \},
\\
C_3 &=& \{ \bar x \mid d(\bar x, C_1) = 1 \} \setminus C_2,
\\  \label{eq:C4}
C_4 &=& V(H^n) \setminus (C_1 \cup C_2 \cup C_3);
\end{eqnarray}
$$ A^j_l (\bar x) = |\{ \bar y \in C_j \mid d(\bar x, \bar y) = l \}|,
\qquad j\in \{1,2,3,4\}, \ \bar x \in V(H^n) $$
(the tuple $(A^i_0 (\bar x),A^i_1 (\bar x),\ldots ,A^i_n (\bar x))$ is known
as the weight distribution of $C_i$ with respect to $\bar x$),

\def\A{\overline{A}\vphantom{A}}
$$ \A^{ij}_l = \frac 1{| C_i|} \sum_{\bar x \in C_i} A^j_l(\bar x),
\qquad i,j\in \{1,2,3,4\}, \ l\in \{0,\ldots,n\}$$
(the tuple $(\A^{ii}_0,\A^{ii}_1,\ldots ,\A^{ii}_n)$ is known
as the inner distance distribution of $C_i$).

Best and Brouwer \cite{BesBro77} showed that $(2^m-3, 2^{n-m}, 3)$ codes are optimal,
i.e. any $(2^m-3, M, 3)$ code satisfies $M\leq 2^{n-m}$.
Moreover,

\begin{lemma}[\cite{BesBro77}]
The inner distance distribution $(\A^{11}_l)_{l=0}^n$
does not depend on the choice of the $(2^m-3, 2^{n-m}, 3)$ code $C_1$.
\end{lemma}

We will also need the following fact:
\begin{lemma}\label{th:SZ} Any $1$-perfect or twofold $1$-perfect code $C$ is antipodal;
i.e., in multiset terms, for any $\bar x \in V(H^n)$ the $C$-multiplicities
of $\bar x$ and $\bar x+\bar 1$ coincide.
\end{lemma}
In the case of $1$-perfect codes this is well-known fact,
which follows from the results \cite{Lloyd,ShSl}.
For twofold $1$-perfect codes, the fact has a similar proof.
Alternatively, Lemma~\ref{th:SZ} follows from the fact that the multiplicity function
of the considered code is, up to an additive constant,
an eigenfunction of $H^n$ with the eigenvalue $-1$
and the corresponding eigenspace has a simple basis from antipodal functions.

\section{An element of equitable partition}
\begin{prop}\label{p:pc}
If $C_1$ is a doubly-shortened $1$-perfect code of length $n$, then
$\A^{11}_n=\A^{24}_{1}=\A^{42}_{1}=0$, $\A^{11}_{n-1}=\A^{44}_{1}=1$, and $A^4_1(\bar x)=\A^{34}_1=2$ for any $\bar x\in C_3$.
\end{prop}
\proofr
Let $C=C_1\times\{00\}  \cup C'\times\{01\}  \cup C''\times\{01\}
\cup C'''\times\{11\}$
be a $1$-perfect code.
If $\bar x \in C_1$ (i.e. $\bar x 00 \in C$),
then $\bar x 00 + \bar 1\in C$, i.e., $\bar x + \bar 1 \in C'''$;
so, $C_2 = C'''$ and $\A^{11}_n = 0$.

If a vertex $\bar y$ is at distance at least $2$ from $C_1$, then, by the definition
of a $1$-perfect code, the vertex $\bar y 00$ is at distance $1$ from an element of $C$,
which is either  $\bar y 01$ or $\bar y 10$. So, $\bar y \in C'\cup C''$. Vise versa,
any $\bar y \in C'\cup C''$ is at distance at least $2$ from $C_1$, because
the minimal distance of $C$ is $3$. So, $C_4=C' \cup C''$.

Because of the minimal distance of $C$, the sets $C_2=C'''$ and $C_4=C' \cup C''$ are at distance more than $1$ from each other. This means $\A^{24}_{1}=\A^{42}_{1}=0$.

We state that for any $\bar x$ from $C_1$
there is exactly one vertex of $C_2$ at the distance $1$ from $\bar x$.
Indeed, the vertex $\bar x 11$ from $V(H^{n+2})$ is at the distance $1$ from
exactly one codeword of $C$, which can be only of type $\bar y 11$,
where $\bar y \in C_2$ and $d(\bar x, \bar y) = 1$. Then the vertex $\bar y+\bar1$ is the only $C_1$-vertex at the distance $n-1$ from $\bar x$; so, $\A^{11}_{n-1} = 1$. The remaining part of the proposition is proved by similar arguments.
\proofend

We will first prove that
\begin{lemma}\label{l:1}
All the numbers
$ \A^{ij}_l$
$( i,j\in \{1,2,3,4\},$ $l\in \{0,\ldots,n\})$
do not depend on the choice of the $(2^m-3, 2^{n-m}, 3)$ code $C_1$.
\end{lemma}
\proofr
Once we have proved that $ \A^{ij}_l$ does not depend on the choice $C_1$,
we know that it is the same as if $C_1$ would be a double-shortened $1$-perfect
(for example, Hamming) code. Moreover if it is equal to the minimal or maximal
possible value of $A^{j}_l(\bar x), \ \bar x \in C_i$, then
$A^{j}_l(\bar x)=\A^{ij}_l$ for any $\bar x \in C_i$.

In particular, for any $\bar x \in  C_1$
$$
A^{1}_n(\bar x)=0 \qquad \mbox{and} \qquad
A^{1}_{n-1}(\bar x)=1.
$$
This means that the sets $C_1$ and $C_2$ are disjoint and
\begin{equation}\label{eq:C1C2}
\mbox{any vertex from
$C_2$ has exactly one neighbor from $C_1$, and vise versa}
\end{equation}
(the fact (\ref{eq:C1C2}) will be used later).
So, $\{C_1,C_2,C_3,C_4\}$ is a partition of $V(H^n)$, and
we can derive relations between the cardinalities of $C_i$:

$$
| C_2| = |C_1|, \quad
| C_3| = (n-1) |C_1|, \quad
| C_4| = |V(H^n)| - |C_1| - |C_2| - |C_3| = 2 |C_1|.
$$

Now we claim the following:
\begin{eqnarray} \label{e:2}
\A^{i2}_l &=& \A^{i1}_{n-l}
\\ \label{e:3}
\A^{i3}_l &=& (n-l+1)\cdot \A^{i1}_{l-l} + (l+1)\cdot \A^{i1}_{l+l} - \A^{i3}_l
\\ \label{e:4}
\A^{i4}_l &=& {n \choose l} - \A^{i1}_l - \A^{i2}_l - \A^{i3}_l
\\ \label{e:1}
|C_i| \cdot \A^{ij}_l &=& |C_j| \cdot \A^{ji}_l
\end{eqnarray}

Indeed, (\ref{e:2}) follows from
$A^{2}_l(\bar x) = A^{1}_{n-l}(\bar x)$,
which is straightforward from the definition of $C_2$;
(\ref{e:3}) follows from
$
A^{3}_l(\bar x) = (n-l+1)\cdot A^{1}_{l-l}(\bar x) + (l+1)\cdot A^{1}_{l+l}(\bar x) - A^{3}_l(\bar x),
$
which is straightforward from the definition of $C_3$ and (\ref{eq:C1C2});
(\ref{e:4}) follows from
\begin{equation}\label{e:1234}
 A^{4}_l(\bar x) + A^{1}_l(\bar x) + A^{2}_l(\bar x) + A^{3}_l(\bar x)={n \choose l},
 \end{equation}
which is from
the fact that $\{C_1,C_1,C_1,C_1\}$ is a partition of $V(H^n)$; the right and left part of (\ref{e:4})
are just different ways to calculate the cardinality of
$\{(\bar x, \bar y) \,|\, \bar x \in C_i, \bar y \in C_j, d(\bar x,\bar y)=l\}$.

Starting from $\A^{11}_l$, we can calculate $\A^{1j}_l$ by (\ref{e:2}-\ref{e:4}),
the values of $\A^{j1}_l$ by (\ref{e:1}),
the values of $\A^{ji}_l$ by (\ref{e:2}-\ref{e:4}); so, Lemma~\ref{l:1} is proved.
\proofend

\begin{theorem}\label{th:pc}
The partition $\{C_1,C_2,C_3,C_4\}$ of $V(H^n)$
is equitable with parameters $$(S_{ij})_{i,j=1}^4 = \left(\begin{array}{cccc}
0&1&n-1&0\\1&0&n-1&0\\1&1&n-4&2\\0&0&n-1&1
\end{array}\right)
.$$
\end{theorem}

\proofr Assume $i,j\in\{1,2,3,4\}$, $\bar x \in C_i$.
We will show that $A^{j}_1(\bar x) = S_{ij}$.

We have already found (\ref{eq:C1C2}) that
$A^{j}_1(\bar x) = 1$ if $(i,j)\in\{(1,2),(2,1)\}$.
Since $C_1$ and $C_2$ are distance-$3$ codes,
$A^{j}_1(\bar x) = 0$ if $(i,j)\in\{(1,1),(2,2)\}$.
Then, by the definition of $C_4$, we have
$A^{j}_1(\bar x) = 0$ if $(i,j)\in \{(1,4),(4,1)\}$.

By (\ref{e:1234}) we get $A^{j}_1(\bar x) = n-0-1-0 = n-1$ if $(i,j)=(1,3)$.

Since $\A^{24}_1=\A^{42}_1=0$, we also have
$A^{j}_1(\bar x) = 0$ if $(i,j)\in \{(2,4),(4,2)\}$.

By (\ref{e:1234}), $A^{j}_1(\bar x) = n-1-0-0 = n-1$ if $(i,j)=(2,3)$.

Let us check that
$A^{j}_1(\bar x) = 1$ if $(i,j)=(4,4)$.
Since $\A^{44}_1 = 1$
(Proposition~\ref{p:pc}),
it is enough to prove that
$A^{j}_1(\bar x)$ is odd.
Indeed, as follows from the arguments above,
the neighborhood of $\bar x$
consists of only $C_3$- and $C_4$-vertices.
Every such $C_3$-vertex is adjacent with exactly
one $C_1$-vertex, which is at distance $2$ from $\bar x$.
While every such $C_1$-vertex is adjacent with exactly
two vertices from the neighborhood of $\bar x$.
So, this neighborhood contains an even number
of vertices from $C_3$ and, consequently,
an odd, from $C_4$.

Automatically,
we get $A^{j}_1(\bar x) = n-0-0-1 = n-1$ if $(i,j)=(4,3)$.

Let us show that
$A^{j}_1(\bar x) = 2$ if $(i,j)=(3,4)$.
%Since $\A^{34}_1 = 2$, it is enough to show that
%$A^{4}_1(\bar x)$ is not greater than $2$.
%
We will calculate the number $T$ of triples
$\{\bar a, \bar b, \bar c\}$ such that $\bar b \in C_3$ is adjacent to both $\bar a,\bar c\in C_4$.
At first, we observe that $T = |C_4|\A^{44}_2$ is independent on the choice of $C_1$.
At second, it can be calculated as
$$\sum_{\bar b\in C_3} \frac{A^4_1(\bar b)(A^4_1(\bar b)-1)}{2} =
\frac12\sum_{\bar b\in C_3} \left(A^4_1(\bar b)\right)^2 - \frac12|C_3|\A^{34}_1;$$
so, by the Cauchy--Bunyakovsky inequality,
$$ T \geq \frac1{2|C_3|} \Big(\sum_{\bar b\in C_3} A^4_1(\bar b)\Big)^2 - \frac{|C_3|}2\A^{34}_1 =
\frac{|C_3|}2  \A^{34}_1(\A^{34}_1-1),
$$
where the equality holds if and only if all $A^4_1(\bar b)$, $\bar b\in C_3$,
are equal to the same value (i.e., to $\A^{34}_1 = 2$). But the last is true when $C_1$
is a doubly-shortened $1$-perfect code
(Proposition~\ref{p:pc}); consequently, it is true for any $(2^m-3, 2^{n-m}, 3)$ code.

%(1) We first check that it is even. Consider
%the neighborhood of $\bar x$. It contains
%one $C_1$-vertex and, by the same arguments as
%above, an even number of $C_2$- and $C_3$-vertices
%(twice the number of $C_1$-vertices at distance $2$ from $\bar x$).
%A remaining even number of vertices are from $C_4$.
%(2) Consider $\bar x \in C_3$ such that
%$A^{4}_1(\bar x)>0$ and adjacent $\bar y \in C_4$.
%There are $\A^{44}_2=(n-1)/2$ vertices of $C_4$
%at the distance $2$ from $\bar y$.
%On the other hand, this number equals
%$\frac 12 \sum_{\bar z\in C_3, d(\bar z,\bar y)=1}(A^{4}_1(\bar z)-1)$ .
%Since $A^{4}_1(\bar z)$ is even and greater than $0$, this value coincides with $(n-1)/2$
%only if $A^{4}_1(\bar z)=2$ for all considered $\bar z$, in particular, for $\bar z=\bar x$.
%So, we have proved that $A^{4}_1(\bar x)\leq 2$ for all $\bar x \in C_3$.

Finally, if $(i,j)=(3,3)$, then $A^{j}_1(\bar x) = n-1-1-2 = n-4$.
\proofend

\begin{note}
 1) If we unify the two parts $C_1$ and $C_2$, say $C_{12} = C_1 \cup C_2$, then
  we will obtain an equitable partition $\{C_{12},C_3,C_4\}$ with parameters
\begin{equation}\label{eq:pc3}
   \left( \begin{array}{ccc} 1 & n-1 & 0 \\ 2 & n-4 & 2 \\ 0 & n-1 & 1
  \end{array} \right).
\end{equation}
  We see that the parameter matrix is symmetrical with respect to interchanging
  of the parts $C_{12}$ and $C_4$. But $C_{12}$ is known to be splittable, while
  the splittability of $C_4$ is questionable. When $C_1$ is a doubly-shortened
  $1$-perfect code, we know that both $C_{12}$ and $C_4$ are splittable.
  Moreover, one can construct an equitable partition with parameters (\ref{eq:pc3})
  whose first and third parts are unsplittable. The problem is if there exists
  such a partition with exactly one of $C_{12}$ and $C_4$ being splittable.

  2) If $C_4$ is splittable, then after splitting it, from the partition $\{C_1,C_2,C_3,C_4\}$
  we obtain an equitable partition with parameters
 $$ \left(\begin{array}{ccccc}
0&1&n-1&0&0\\1&0&n-1&0&0\\1&1&n-4&1&1\\0&0&n-1&0&1\\0&0&n-1&1&0
\end{array}\right), $$
which also have some obvious symmetries.
\end{note}
\begin{note}
An equitable partition with $\{C_{12},C_3,C_4\}$ of $H^n$ with parameters (\ref{eq:pc3})
generates an equitable partition $\{G_1,G_2,G_3,G_4\}$ of $H^{n'}$, $n'=n+1$
with parameters
 $$ \left(\begin{array}{cccc}
0&n'&0&0\\2&0&n'-2&0\\0&n'-2&0&2\\0&0&n'&0
\end{array}\right) $$
 as follows:
 \begin{eqnarray*}
   G_1 &=& \{ \bar x \alpha \mid \bar x=x_1x_2...x_n\in C_{12},\ \alpha=x_1+...+x_n\bmod 2 \}, \\
   G_4 &=& \{ \bar x \beta \mid \bar x=x_1x_2...x_n\in C_{4},\ \beta=x_1+...+x_n+1\bmod 2 \},    \\
   G_2 &=& \{ \bar y \in V(H^{n'}) \mid d(\bar y, C_1)=1 \}, \\
   G_3 &=& \{ \bar y \in V(H^{n'}) \mid d(\bar y, C_4)=1 \}.
 \end{eqnarray*}
 This partition can be viewed as an ``extended'' version of the partition $\{C_{12},C_3,C_4\}$; the spittability of $C_{12}$ or $C_4$ is equivalent to the spittability of $G_{1}$ or $G_4$ respectively. But the distance $1$ between vertices of, say, $C_4$ corresponds to the
 distance $2$ between the corresponding vertices of $G_4$; and the graph of distances $1$ and $2$ of $C_4$ corresponds to the graph of distances $2$ of $G_4$, which emphasize the ``equal status'' of the all edges of the graph.
 Of cause if $G_1$ and/or $G_4$ are splittable, then splitting gives an equitable partition
 of $H^{n'}$ into $6$\,/\,$5$ parts with corresponding parameters.
 If both $G_1$ and $G_4$ are splittable (say, into $G'_1$, $G''_1$ and $G'_4$, $G''_4$ respectively), then the equitable partition $\{G_1,G_2,G_3,G_4\}$
 (defined in some other terms) is also known as a
 \emph{code-generating factorization} of $H^{n'}$ \cite{VasSol97}. Indeed, the code $G'_10\cup G''_11\cup G'_40\cup G''_41$ is $1$-perfect.
\end{note}

\section{Embedding in twofold $1$-perfect codes}\label{s:21p}
\begin{theorem}\label{th:twofold2}
Let $C_1$ be a $(n=2^m-3, 2^{2^m-m-3},3)$ code.
Then the set $ C_1 00 = \{\bar x 00 \mid \bar x \in C_1 \}$
is a subset of a unique twofold $1$-perfect
code $B$ with the following properties:\\
{\rm a)} the multiplicity of any codeword of $ C_1 00$ is $2$;\\
{\rm b)} any codeword $\bar x$ with the last two symbols $01$ or $10$ satisfies $\bar x + 0...011 \in C$.
\end{theorem}
\proofr
{\em Existence.}
Let $B=2*C_1 00 \cup 2*C_2 11 \cup C_401 \cup C_4 10$.
Obviously, $B$ satisfies a), b), and $ C_1 00\subset B$.
The fact that $B$ is a twofold $1$-perfect code
is straightforward from Theorem~\ref{th:pc};
we leave the details as an exercise.

{\em Uniqueness.} Assume $B$ is a twofold $1$-perfect
code satisfying a), b), and $ C_1 00\subset B$.
Define
\begin{eqnarray*}
C_2= \{ \bar x \mid \bar x11 \in B\},
\\
C_4= \{ \bar x \mid \bar x01 \in B\},
\\
C_3=V(H^n)\setminus(C_1\cup C_2\cup C_4)
\end{eqnarray*}
From the antipodality of $B$, we have $C_2=C_1$.
As follows from the definition of twofold $1$-perfect
codes, any codeword of multiplicity $2$ cannot be at distance
$1$ or $2$ from any other codeword. Consequently,
1) the distance between $C_1$ and $C_4$, as well as
between $C_2$ and $C_4$, cannot be less than $2$;
2) the multiplicity of the words of form $\bar x01$ in $B$
is less than $2$.

Now we see that, by numerical reasons, $C_4$ consists
of the all vertices at the distance more than $1$ from
$C_1$. Thus, $C_2$, $C_3$, and $C_4$ satisfy
(\ref{eq:C2})-(\ref{eq:C4}), and $B$ is unique.
\proofend

By similar arguments, the following is also true:
\begin{theorem}\label{th:twofold1}
Let $C_1$ be a $(n=2^m-3, 2^{2^m-m-3},3)$ code.
Then the set $C_1 00$
is a subset of a unique twofold $1$-perfect
code $D$ whose all codewords $\bar x$ satisfy
$\bar x + 0...011 \in C$.
\end{theorem}

\section{Embedding in $1$-perfect codes}\label{s:1p}
\begin{theorem}\label{th:1perf}
Let $C_1$ be a $(n=2^m-3, 2^{2^m-m-3},3)$ code.
The following four statements are mutually equivalent:\\
{\rm a)} the set $C_1 00$ is a subset of a $1$-perfect code $C$;\\
{\rm b)} the set $C_4$ defined in {\rm (\ref{eq:C4})} is splittable;\\
{\rm c)} the twofold $1$-perfect $B$ from Theorem~{\rm\ref{th:twofold2}} is splittable;\\
{\rm d)} the twofold $1$-perfect $D$ from Theorem~{\rm\ref{th:twofold1}} is splittable.
\end{theorem}
\proofr
Clearly, each of c) and d) implies a).

Since $B$, as well as $D$, includes $C_4 01$, each of c) and d) implies b).

Conversely, assume b) holds and $C_4 = C' \cup C''$ where $C'$ and $C''$ are distance-$3$ codes. Then
\begin{eqnarray*}
B & = & (C_1 00 \cup C_4 11 \cup C'01 \cup C''10)\cup (C_1 00 \cup C_4 11 \cup C'10 \cup C''01) , \\
D& = & (C_1 00 \cup C_4 11 \cup C'01 \cup C''10)\cup (C_4 00 \cup C_1 11 \cup C'10 \cup C''01),
\end{eqnarray*}
and c), d) hold.

Assume a) is true. Define $C' = \{\bar x \mid \bar x 01\in C\}$ and $C'' = \{\bar x \mid \bar x 10\in C\}$. Because of the code distance $3$ of $C$, we see that $C'$ and $C''$ are disjoint
and at the distance at least $2$ from $C_1$. So, since $|C'|+|C''|=|C_4|$,
we get $C_4=C'\cup C''$, and b) holds.
\proofend

\begin{note} The splittability of any of the sets $C_4$, $B$, $D$ is equivalent to the biparticity
of its graph of distances $1$ and $2$ (two codewords $\bar x$ and $\bar y$
 are adjacent if and only if $d(\bar x,\bar y)\in\{1,2\}$).
 In this graph for $D$, the vertices of types $\bar x00$ and $\bar x11$
 are not connected with the vertices of types $\bar x01$ and $\bar x10$,
 and the subgraph generated by the former vertices is bipartite,
 while the biparticity of the remaining subgraph is questionable.
 In $B$, the codewords of types $\bar x00$ and $\bar x11$ have the multiplicity $2$,
 and they are isolated in the graph of distances $1$ and $2$.
\end{note}

\begin{note}
If $\nu$ is the number of connected components in the graph of distances $1$ and $2$ of $C_4$,
then the number of different $1$-perfect codes including $C_1 00$ is $2^{\nu}$.
As follows from the tight lower bound on the size of the difference between two $1$-perfect codes \cite{Sol90en,EV:94}, the cardinality of a connected component is not less than $2^{\frac{n-1}2}$,
and so $\nu \geq \frac{2^{\frac{n-3}2}}{n+1}$. If $C_1$ is linear, then $\nu$ achieves this bound.
\end{note}

\section{Unsplittable twofold STS}\label{s:2STS}
If we consider a $1$-perfect code containing $\bar 0$, then all the weight-$3$ codewords
compose a design known as a Steiner triple system, or STS.
The characteristic property of an \emph{STS} is that every weight-$2$ word is at distance $1$
from exactly one word of the STS. (Strictly speaking, an STS is defined as a pair $(V,B)$,
where $V$ is some set and $B$ is a collection of $3$-subsets of $V$, named blocks,
such that every $2$-subset of $V$ is included in exactly one block.)

If we consider a twofold $1$-perfect code $C$ such that the multiplicity of $\bar 0$ is $2$,
 then all the weight-$3$ codewords compose a design, which can be called a twofold STS.
 If $C$ comes from Theorem~\ref{th:twofold2}, then the corresponding STS satisfies

 a) any codeword of type $\bar x 00$ or $\bar x 11$ has the multiplicity $2$;

 b) for any $\bar x$ of the corresponding length, $\bar x 01$ and $\bar x 10$
 are codewords or not simultaneously.

For the length $15$, there exists a twofold STS meeting a) and b) that cannot be split
 into two STS. This fact has not direct connection with the problem considered in this paper:
 on one hand, it is not proved that there exists a twofold $1$-perfect code that include this STS
 (e.g., for the length $15$,
 there exist STSs that are not embeddable in a $1$-perfect code \cite{OstPot2007});
 on the other hand, the splittability of the all twofold STS included in a twofold $1$-perfect
 code would not mean the splittability of the twofold $1$-perfect itself. Nevertheless,
 the existence of such an object seems to be interesting. The following is the list of the words
 of the mentioned example (the unsplittability follows from the existence of a $5$-cycle
 in the distance-$2$ graph):
\texttt{\small\mbox{}\\
\quad0000000\,00000\,1\,11$\,\times 2$\,, \quad \\[-0.5ex]
\quad0000000\,11000\,0\,01\,, \quad %\\[-1.0ex]
\quad0000000\,11000\,0\,10\,, \quad \\[-1.0ex]
\quad0000000\,01100\,0\,01\,, \quad %\\[-1.0ex]
\quad0000000\,01100\,0\,10\,, \quad \\[-1.0ex]
\quad0000000\,00110\,0\,01\,, \quad %\\[-1.0ex]
\quad0000000\,00110\,0\,10\,, \quad \\[-1.0ex]
\quad0000000\,00011\,0\,01\,, \quad %\\[-1.0ex]
\quad0000000\,00011\,0\,10\,, \quad \\[-1.0ex]
\quad0000000\,10001\,0\,01\,, \quad %\\[-1.0ex]
\quad0000000\,10001\,0\,10\,, \quad \\[-1.0ex]
\quad1100000\,00000\,0\,01\,, \quad %\\[-1.0ex]
\quad1100000\,00000\,0\,10\,, \quad \\[-1.0ex]
\quad0110000\,00000\,0\,01\,, \quad %\\[-1.0ex]
\quad0110000\,00000\,0\,10\,, \quad \\[-1.0ex]
\quad0011000\,00000\,0\,01\,, \quad %\\[-1.0ex]
\quad0011000\,00000\,0\,10\,, \quad \\[-1.0ex]
\quad0001100\,00000\,0\,01\,, \quad %\\[-1.0ex]
\quad0001100\,00000\,0\,10\,, \quad \\[-1.0ex]
\quad0000110\,00000\,0\,01\,, \quad %\\[-1.0ex]
\quad0000110\,00000\,0\,10\,, \quad \\[-1.0ex]
\quad0000011\,00000\,0\,01\,, \quad %\\[-1.0ex]
\quad0000011\,00000\,0\,10\,, \quad \\[-1.0ex]
\quad1000001\,00000\,0\,01\,, \quad %\\[-1.0ex]
\quad1000001\,00000\,0\,10\,, \quad \\[-0.5ex]
\quad0000000\,01010\,1\,00$\,\times 2$\,, \quad \\[-1.0ex]
\quad0100010\,00000\,1\,00$\,\times 2$\,, \quad \\[-1.0ex]
\quad0010100\,00000\,1\,00$\,\times 2$\,, \quad \\[-1.0ex]
\quad0001000\,00100\,1\,00$\,\times 2$\,, \quad \\[-1.0ex]
\quad1000000\,10000\,1\,00$\,\times 2$\,, \quad %\\[-1.0ex]
\quad0000001\,00001\,1\,00$\,\times 2$\,, \quad \\[-1.0ex]
\quad0100000\,10100\,0\,00$\,\times 2$\,, \quad %\\[-1.0ex]
\quad0000010\,00101\,0\,00$\,\times 2$\,, \quad \\[-1.0ex]
\quad0010000\,10010\,0\,00$\,\times 2$\,, \quad %\\[-1.0ex]
\quad0000100\,01001\,0\,00$\,\times 2$\,, \quad \\[-1.0ex]
\quad1010000\,00001\,0\,00$\,\times 2$\,, \quad %\\[-1.0ex]
\quad0000101\,10000\,0\,00$\,\times 2$\,, \quad \\[-1.0ex]
\quad0101000\,00001\,0\,00$\,\times 2$\,, \quad %\\[-1.0ex]
\quad0001010\,10000\,0\,00$\,\times 2$\,, \quad \\[-1.0ex]
\quad1000010\,00010\,0\,00$\,\times 2$\,, \quad %\\[-1.0ex]
\quad0100001\,01000\,0\,00$\,\times 2$\,, \quad \\[-1.0ex]
\quad1001000\,01000\,0\,00$\,\times 2$\,, \quad %\\[-1.0ex]
\quad0001001\,00010\,0\,00$\,\times 2$\,, \quad \\[-1.0ex]
\quad0100100\,00010\,0\,00$\,\times 2$\,, \quad %\\[-1.0ex]
\quad0010010\,01000\,0\,00$\,\times 2$\,, \quad \\[-1.0ex]
\quad1000100\,00100\,0\,00$\,\times 2$\,, \quad %\\[-1.0ex]
\quad0010001\,00100\,0\,00$\,\times 2$.
}

\section{MDS codes and double-MDS-codes}\label{s:2MDS}
Let $Q^m=(V(Q^m),E(Q^m))$ denotes the graph whose vertex set is the set $\{0,1,2,3\}^m$ of quaternary $n$-words,
two words being adjacent if and only if they differ in exactly one position.
By a \emph{$4$-clique} we mean a set of four words of $V(Q^m)$ differing in exactly one position.

A subset $M$ of $V(Q^m)$ is called an MDS code (with distance $2$) if every
\emph{$4$-clique} contains exactly one word of $M$. Equivalently, $M$ is a
distance $2$ code of cardinality $4^{m-1}$. Equivalently, $\{M,V(Q^m)\setminus M\}$
is an equitable partition of $Q^m$ with parameter matrix $((0,3m)(m,2m))$.

We call a subset $M$ of $V(Q^m)$ a double-MDS-code if every
\emph{$4$-clique} contains exactly two word of $M$.
Equivalently, $\{M,V(Q^m)\setminus M\}$
is an equitable partition of $Q^m$ with parameter matrix $((m,2m)(2m,m))$.

A double-MDS-code is \emph{splittable} if it is the union of two (disjoint) MDS codes.

Denote $P_0 = \{0000,1111\}$, $P_1 = \{0011,1100\}$, $P_2 = \{0101,1010\}$,
$P_3 = \{0110,1001\}$ $\subset V(H^4)$ and
$P'_0 = \{000,111\}$, $P'_1 = \{011,100\}$, $P'_2 = \{101,010\}$,
$P'_3 = \{110,001\}$ $\subset V(H^3)$.
Let $C\in V(H^{m-1})$ be a $1$-perfect binary code; denote
$C^* = \{000c_1\,000c_2...000c_{m-1}\,000 \mid c_1c_2...c_{m-1}\in C\}$.
For any subset $M$ of $V(Q^m)$ we define the code $S(M)\subset V(H^{4m-1})$ as follows:
\begin{equation}\label{eq:GCC}
    S(M) = \bigcup_{\mu_1...\mu_{m} \in M}P_{\mu_1}P_{\mu_2}...P_{\mu_{m-1}}P'_{\mu_{m}}
    + C^*
\end{equation}
Here, for two sets of words $P^1\subset V(H^r)$ and $P^2\subset V(H^l)$, $P^1P^2 = \{x_1...x_ry_1...y_l \mid x_1...x_r\in P^1,\ y_1...y_l\in P^2\}$;
and if $r=l$, then $P^1+P^2 = \{z_1...z_r \mid z_i=x_i+y_i\bmod 2,\  x_1...x_r\in P^1,\ y_1...y_r\in P^2\}$.

\begin{prop}\label{p:mds-perf}
1)  If the code distance of $M$ is not less than $2$,
then the code distance of $S(M)$ is at least $3$;
if $M$ is an MDS code, then $S(M)$ is a $1$-perfect code.
2)  If $M$ is a splittable (unsplittable) double-MDS-code, then $S(M)$ is a splittable (unsplittable) twofold $1$-perfect code.
\end{prop}
\proofr[(a sketch)]
P. 1) is proved in \cite{Phelps84}, in more general form.

Similarly, if $M$ is a double-MDS-code, then $S(M)$ is a twofold $1$-perfect code
(it is straightforward to check the definition). If $M=M'\cup M''$ for some MDS codes
$M'$ and $M''$, then $S(M)=S(M')\cup S(M'')$, where $S(M')$ and $S(M'')$ are $1$-perfect codes. Otherwise, the distance-$1$ graph of $M$ has an odd cycle, and it is easy to find a corresponding cycle of the same length in the graph of distances $1$ and $2$ of $S(M)$, which implies that $S(M)$ is unsplittable.
\proofend

\begin{theorem}
  Let $m=2^{k-2}$. Assume there exists an unsplittable double-MDS-code $M_1\subset V(Q^{m-1})$ such that the double-MDS-code $M_0=V(Q^{m-1})\setminus M_1 $ is splittable.
  Then there exist a $(n=2^k-3, 2^{2^k-k-3},3)$ code $C_1$ such that $C_100$
  is not a subset of a $1$-perfect code.
\end{theorem}
\proofr
Let $M_0 = M' \cup M''$ where $M'$ and $M''$ are disjoint MDS codes.

Denote $M=M_0 0 \cup M_0 1 \cup M_1 2 \cup M_1 3\subset V(Q^{m})$.
By the definition, $M$ is a double-MDS-code.
Since $M_1$ is unsplittable, $M$ is unsplittable too.
Then, by Proposition~\ref{p:mds-perf}, the set
$$D = S(M)$$ is an unsplittable twofold $1$-perfect code.

Now, consider the set
$$C = S(M'0\cup M'1).$$
Since the code distance of $M'0\cup M'1$ is $2$,
the code distance of $C$ is at least $3$,
by Proposition~\ref{p:mds-perf}.
Half of the codewords of $C$ have $00$ in the last two positions (the others, $11$);
let $C_100$ denote the corresponding subcode.

We have: $|C_1|=\frac18|C|=2^{2^k-k-3}$; the code distance of $C_1$ is $3$;
$C_100\subset D$ where $D$ is an unsplittable twofold $1$-perfect code
whose all codewords $\bar x$ satisfy
$\bar x + 0...011 \in C$.
By Theorems~\ref{th:twofold1} and~\ref{th:1perf}, the proof is over.
\proofend

\textbf{Conjecture} (V. Potapov). \emph{Any double-MDS-code $M$ in $Q^m$ is splittable if and only if its complement $V(Q^m)\setminus M$ is splittable.}

This is equivalent to the following statement:
\emph{Any  $4\times 4 \times \ldots 4 \times 2$ latin hypercuboid
is completable to a  $4\times 4 \times \ldots 4 \times 4$ latin hypercube.}
A $q\times q \times \ldots q \times p$ \emph{latin hypercuboid} of order $q$ (if $p=q$, \emph{latin hypercube})
is a function $f:\{0,\ldots,q-1\}^{m-1}\times\{0,\ldots,p-1\}\to \{0,\ldots,q-1\}$ such that
$f(\bar x)\neq f(\bar y)$ for any $\bar x$ and $\bar y$ differing in exactly one position.
Examples of non-completable latin cuboids are constructed in \cite{Kochol,MK-W:small}

Another equivalent formulation: \emph{Let $K_4^m$ be the direct product of $m$ copies of the complete
graph on $4$ vertices. If $V(K_4^m)$ is partitioned into two subsets that generate subgraphs of degree $m$, then these subgraphs are bipartite or not bipartite simultaneously.}

It seems perspective to use the characterization of the distance-$2$ MDS codes over the quaternary alphabet (latin hypercubes of order $4$) \cite{KroPot:4} to prove this conjecture. 
Nevertheless, the analysis of all subcases needs some work, which is not completed at this moment. In any case, it is interesting to find an independent proof.

\nocite{KroPot:nonsplittable}
%\bibliographystyle{plain}
%\bibliography{../k}

\providecommand\href[2]{#2} \providecommand\url[1]{\href{#1}{#1}} \def\DOI#1{
  {DOI}: \href{http://dx.doi.org/#1}{#1}}

\end{document}